\documentstyle[a4,leqno,amsfonts,12pt]{article}
\newcommand{\C}{{\bf C}}

\newcommand{\N}{{\bf N}}

\newcommand{\D}{{\bf D}}
\newcommand{\de}{\delta}

\newcommand{\mb}{\mbox}

\newcommand{\beq}{\begin{equation}}
\newcommand{\eeq}{\end{equation}}
\newcommand{\oge}{\succeq}
\newcommand{\ole}{\preceq}
\newcommand{\ve}{\varepsilon}

\newcommand{\ov}{\overline}
\newcommand{\al}{\alpha}

\newcommand{\Om}{\Omega}
\newcommand{\om}{\omega}
\newcommand{\z}{\zeta}

\newcommand{\kap}{\mb{ cap}}
\newcommand{\ga}{\gamma}

\newtheorem{th}{Theorem}

\newtheorem{lem}{Lemma}

\newcommand{\ueberschrift}{\bigskip\goodbreak\noindent\bigskip}
\newcounter{theabsatz}
\newcommand{\absatz}[1]{\stepcounter{theabsatz} \ueberschrift
                           {\large \bf \arabic{theabsatz}. {#1}} \setcounter{equation}{0}}

\begin{document}
\mathsurround=2pt

\begin{center}
{\large \bf CHEBYSHEV POLYNOMIALS on a SYSTEM
of CONTINUA}\\[3ex]

 Vladimir V.
Andrievskii \\[3ex]

\end{center}

\begin{abstract}
 The estimates of the uniform norm of the Chebyshev polynomial
 associated with a compact set $K$ consisting of a finite
 number of continua in the complex plane are established.
 These estimates are exact (up to a constant factor)
 in the case where the components of $K$ are either quasismooth
 (in the sense of Lavrentiev) arcs or closed Jordan domains
 bounded by a quasismooth curve.
\end{abstract}

\footnotetext{
Date received:$\quad$. Communicated by


 }

\absatz{Introduction and Main Results}

Let $K\subset\C$ be a compact set in the complex plane
consisting of disjoint  closed connected sets (continua) $K^j,j=1,2,\ldots,m$,
i.e.,
$$
K=\cup_{j=1}^mK^j;\quad K^j\cap K^k=\emptyset\mb{ for }j\neq k;
\quad\mb{ diam}(K^j)>0.
$$
Here
$$
\mb{diam}(S):=\sup_{z,\z\in S}|z-\z|,\quad S\subset\C.
$$
Denote by $T_n(z)=T_n(z,K),n\in\N:=\{ 1,2,\ldots\}$
the $n$-th Chebyshev polynomial associated with $K$, i.e.,
$T_n(z)=z^n+c_{n-1}z^{n-1}+\ldots+c_0, c_k\in\C$ is the (unique)
monic polynomial which minimizes the supremum norm $||T_n||_K:=
\sup_{z\in K}|T_n(z)|$ among all monic polynomials of the same
degree.

It is well-known (see, for example, \cite[Theorem 5.5.4 and
Corollary 5.5.5]{ran}) that
$$
||T_n||_K\ge\kap(K)^n\quad\mb{and}\quad
\lim_{n\to\infty}||T_n||_K^{1/n}=\kap(K),
$$
where$\kap(S)$ denotes the logarithmic capacity of a compact set
$S\subset\C$ (see \cite{ran, saftot}).

We are interested in finding the estimates of $||T_n||_K$ from above.
This problem attracted attention of many mathematicians;
for the complete survey of results and further citations, see
\cite{smileb}-\cite{wid}.
Without loss of generality we always  assume that $\Om:=\C\setminus K$ is connected.
Denote by $c, c_1,\ldots$ positive constants
(different in different sections) that are either absolute or they depend only
on $K$; otherwise, the dependence on  other parameters is explicitly
stated.
\begin{th}\label{th1}
Under the above assumptions,
\beq\label{1.1}
||T_n||_K\le c_1\log(n+1){\em{\kap}}(K)^n,\quad n\in\N.
\eeq\end{th}
 The estimate (\ref{1.1}) is surprising. It is unexpected even in the case where
 $K$ is a continuum, i.e., $m=1$. In this case, one of the major sources
 for estimates of  $||T_n||_K$ are Faber polynomials $F_n(z)=F_n(z,K)$
 associated with $K$ (see \cite{smileb, sue}). According to Koevary
 and Pommerenke \cite[Theorem 1]{kovpom} (see also \cite[Chapter IX, \S 4]{sue}),
 there exist a continuum $K_0$ with$\kap(K_0)=1$ and an infinite set
 $\Lambda\subset\N$ such that for the (monic) polynomial $F_n(z)=F_n(z,K_0)$
 and $n\in\Lambda$ we have
 $$
 ||F_n||_{K_0}> n^\al,\quad \al=0.138.
 $$
It is interesting to compare this classic estimate with (\ref{1.1}).

If more information is known about the geometry of $K$, (\ref{1.1}) can be improved in the following
way.
A Jordan curve $L\subset\C$ is called quasismooth
(in the sense of
Lavrentiev) (see \cite[p. 163]{pom}) if
for every $z_1,z_2\in L$,
\beq\label{1.3}
|L(z_1,z_2)|\le c_2|z_2-z_1|,
\eeq
where $L(z_1,z_2)$ is the shorter arc of $L$ between $z_1$ and $z_2$,
 a constant $c_2$ depends on $L$,  and $|S|$ is the linear measure (length) of a
(Borel) set $S\subset\C$ (see \cite[p. 129]{pom}). Any subarc of a quasismooth
curve is called a quasismooth arc.
\begin{th}\label{th2} Let each $K^j$ in the definition of $K$
be either  a quasismooth arc
or a closed Jordan domain bounded by  a quasismooth curve. Then
\beq\label{1.4}
||T_n||_K\le c_3 {\em{\kap}}(K)^n,\quad n\in\N.
\eeq
\end{th}
In the case of sufficiently smooth components $K^j$, the estimate (\ref{1.4})
was proved in \cite[Theorem 11.4]{wid} and recently in \cite[Theorem 1.3]{tot14}.
Comparing (\ref{1.4}) with known estimates for the Faber polynomials
for a domain with the quasismooth boundary  (see \cite{lesvinwar})
we also see the improvement by a logarithmic factor.

Our constructions below are based on the method of discretization of
the equilibrium measure due to Totik \cite{tot12}-\cite{tot14}, representation
of the Green function via special conformal mappings due to Widom \cite{wid},
and distortion properties of conformal mappings which can be found, for
example, in \cite{pom} or \cite{andbla}.

We use the following notation.
For the functions $a>0$ and $b>0$, we write
$a\ole b$ (order inequality) if $a\le cb$. The expression $a\asymp b$
means that $a\ole b$ and $b\ole a$ simultaneously.
Moreover, for $z\in\C$ and $\de>0$ we let
$$
C(\de):=\{z:|z|=\de\}, \quad \D:=\{ z:|z|<1\},\quad \D^*:=\C\setminus\ov{\D},
$$
$$
d(S_1,S_2):=\inf_{z_1\in S_1,z_2\in S_2}|z_2-z_1|,\quad S_1,S_2\subset\C.
$$
Let $m_2(S)$ be the two-dimensional Lebesgue measure (area) of a (Borel) set $S\subset\C$.
For a bounded Jordan curve $J\subset\C$, denote by  int$(J)$  the bounded
 component of $\C\setminus J$.

\absatz{Totik-Type Polynomials}

In this section we review (in more general setting)
 the construction of  the monic polynomials suggested
in \cite{tot12}-\cite{tot14}. Let $K=\cup _{j=1}^mK^j$ be as above. We recall some general
facts from potential theory which can be found, for example,
in \cite{wal, wid, ran, saftot}.
 Denote  by $g(z)=g_\Om(z,\infty), z\in \Om$  Green's function for
 $\Om$ with pole at infinity. It has a multiple-valued harmonic conjugate
 $\tilde{g}(z)$. Let
 $$
 \Phi(z):=\exp(g(z)+i\tilde{g}(z)),
 $$
 $$
 K_s:=\{z\in\Om:g(z)=s\},\quad s>0.
 $$
 Note that
 \beq\label{2.1}
 \kap(K_s)=e^s\kap(K).
 \eeq
 Let $s_0>0$  be such that for $0<s<2s_0$, the set
 $K_s=\cup_{j=1}^mK_s^j$ consists  of $m$ mutually disjoined
 curves, where $K_s^j$ is the curve surrounding $K^j$.
 Moreover, we can fix $s_0$ so small  that
\beq\label{2.1n}
 d(\z,K^j)=d(\z,K),\quad \z\in \mb{int}(K^j_{2s_0}).
 \eeq
 Let $\mu=\mu_K$ be the equilibrium measure of $K$.
 According to Gauss' Theorem (see \cite[p. 83]{saftot}),
 for the net-change of the function
 $\tilde{g}$, we obtain
 $$
 \Delta_{K_s^j}\tilde{g}:=\int_{K_s^j}\frac{\partial \tilde{g}(\z)}{\partial t_\z}
 |d\z|=2\pi\om_j
 $$
 (see \cite[p. 140]{wid}), where $0<s<2s_0$, $t_\z$ is the tangent vector to the
 curve $K_s^j$ (traversed in the positive, i.e., counterclockwise
  direction) at $\z$, and $\om_j:=\mu(K^j).$

 Therefore, the function
 $
 \phi_j:=\Phi^{1/\om_j}(\z)
 $
is a conformal and univalent mapping of
$
\Om^j:=\mb{int}(K_{2s_0}^j)\setminus K^j
$
 onto the annulus
 $
 A^j:=\{w:1<|w|<e^{2s_0/\om_j}\}
 $
 as well as
 $$
 K_s^j=\{\z\in \Om^j:|\phi_j(\z)|=e^{s/\om_j}\},\quad 0<s<2s_0.
 $$

 Let $\mu_s:=\mu_{K_s} $ be the equilibrium measure of $K_s.$
 By Gauss' Theorem (see \cite[p. 83]{saftot}),
 $$
 \mu_s(K_s^j)=\mu(K^j)=\om_j,\quad 0<s<2s_0.
 $$
 Moreover, by virtue of \cite[p. 90, Theorem 1.4]{saftot}, for any arc
 $$
 \ga=\{\z\in K_s^j:\theta_1\le\arg \phi_j(\z)\le\theta_2\},
 \quad 0<\theta_2-\theta_1\le 2\pi
 $$
 we have
 \beq\label{2.5}
 \mu_s(\ga)=\frac{(\theta_2-\theta_1)\om_j}{2\pi}\, .
 \eeq
 Assuming that $n\in\N$ is sufficiently large, i.e.,
 $n>n_1:=10/(\min_{j}\om_j)$ we let
 $$
 n_j:=[n\om_j],\quad j=1,2,\ldots,m-1,
 $$
 $$
 n_m:=n-(n_1+\ldots+n_{m-1}),
 $$
 where $[a]$ means the integer part of  a real number $a$.

 Therefore,
 \beq\label{2.6}
 0\le n_m-n\om_m=\sum_{j=1}^{m-1}(n\om_j-n_j)\le m-1.
 \eeq
 Next, for $0<s<s_0/2$, we represent each $K^j_s$
 as the union of closed subarcs $I^j_k,k=1,\ldots,n_j$ such that
 their interiors do not intersect,
 $$
 I_k^j\cap I_{k+1}^j=:\xi_k^j,\quad k=1,\ldots, n_j-1,
 $$
 and $I_{n_j}^j\cap I_{1}^j=:\xi_{n_j}^j$ are points of $K_s^j$
 ordered in a positive  direction, and
 $$
 \mu_s(I_k^j)=\frac{\om_j}{n_j},\quad k=1,\ldots,n_j.
 $$
  Let $\xi_0^j:=\xi^j_{n_j}$ and let $\eta_0^j<\eta_1^j<\ldots<\eta_{n_j}^j=\eta_0^j+2\pi$
  be determined by $\phi_j(\xi_k^j)=\exp((s/\om_j)+i\eta_k^j)$, i.e.,
  $\eta_k^j-\eta_{k-1}^j=2\pi/n_j$.
  Consider
  $$
  \tilde{B}_k^j:=\left\{ w=re^{i\eta}:\eta_{k-1}^j\le\eta\le\eta_k^j,
  0\le r-e^{s/\om_j}\le\frac{2\pi}{n_j}\right\}.
  $$
 Further, we assume that
  \beq\label{2.d}
 n>n_2:=\frac{2\pi}{\min_j\left( e^{s_0/\om_j}- e^{s_0/(2\om_j)}\right)},
 \eeq
 which implies
 $
\tilde{B}_k^j\subset E^j:=\{ w:1<|w|<e^{s_0/\om_j}\}.
 $
Let $B_k^j:=\psi_j(\tilde{B}_k^j)$, where $\psi_j:=\phi_j^{-1}$.
\begin{lem}\label{lem2.1}
There exist  constants $c$ and $n_0=n_0(K,c)\in\N$ such that for $s=s(n)=c/n<s_0/2$ and
$n>n_0$ we have
\beq\label{2.7}
d(I_k^j,K)\ole|\xi_{k}^j-\xi_{k-1}^j|\le \mb{\em diam}(I_k^j)\le
|I_k^j|\le 0.1 d(I_k^j,K)
\eeq
\beq\label{2.8}
|\xi_{k}^j-\xi_{k-1}^j|^2\asymp
m_2(B_k^j)\asymp \mb{\em diam}(B_k^j)^2\le \frac{1}{4}d(B_{k}^j,K)^2.
\eeq
\end{lem}
For the proof of the lemma, see Section 3.

For sufficiently large $n>n_0$ and $s:=c/n$ as in Lemma \ref{lem2.1}, consider the points
$$
\z_k^j:=\frac{1}{\mu_s(I_k^j)}\int_{I_k^j}\xi d\mu_s(\xi)
$$
and the polynomial
$$
P_n(z):=\prod_{j=1}^m\prod_{k=1}^{n_j}(z-\z_k^j).
$$
According to Lemma \ref{lem2.1},
$$
|\z_k^j-\xi_k^j|=\left|\frac{1}{\mu_s(I_k^j)}\int_{I_k^j}(\xi-\xi_k^j)d\mu_s(\xi)
\right|
\le 0.1d(I_k^j,K)
$$
and for $\xi\in I_k^j$,
\beq\label{2.333}
|\xi-\z_k^j|\le |\xi-\xi_k^j|+|\xi_k^j-\z_k^j|\le 0.2d(I_k^j,K).
\eeq
For $z\in K$ we have
\begin{eqnarray}
n\log\kap(K_s)&=&\sum_{j=1}^m\sum_{k=1}^{n_j}\left( n-\frac{n_j}{\om_j}\right)
\int_{I_k^j}\log|z-\xi|d\mu_s(\xi)
\nonumber\\
&& +\sum_{j=1}^m\sum_{k=1}^{n_j}  \frac{1}{\mu_s(I_k^j)}
\int_{I_k^j}\log|z-\xi|d\mu_s(\xi)\nonumber\\
\label{2.9}
&=:&\Sigma_1(z)+\Sigma_2(z).
\end{eqnarray}
Since for $\xi\in K_s$ and $z\in K$,
\begin{eqnarray*}
\int_{K_s}|\log|z-\xi||d\mu_s(\xi)&\le&|\log\mb{diam}(K_s)|
+\int_{K_s}\log\frac{\mb{diam}(K_s)}{|z-\xi|} d\mu_s(\xi)\\
&\le &2\log^+\mb{diam}(K_s)-\log\kap(K_s)\ole 1,
\end{eqnarray*}
by virtue of (\ref{2.6}) we obtain
\beq\label{2.10}
|\Sigma_1(z)|\ole \int_{K_s}|\log|z-\xi||d\mu_s(\xi)\ole 1.
\eeq
Furthermore,
\beq\label{2.11}
\log|P_n(z)|-\Sigma_2(z)=\sum_{j=1}^m\sum_{k=1}^{n_j}\frac{1}{\mu_s(I_k^j)}
\int_{I_k^j}\log\left|\frac{z-\z_k^j}{z-\xi}\right|d\mu_s(\xi)
\eeq
which, together with (\ref{2.7}) and (\ref{2.333}), imply that
for $z\in K$ and $\xi\in I_k^j$,
\begin{eqnarray}
\log\left|\frac{z-\z_k^j}{z-\xi}\right|&=&-\Re \log\left(1+\frac{\z_k^j-\xi}{z-\z_k^j}\right)
\nonumber\\
\label{2.12}
&=&\Re \left(\frac{\xi-\z_k^j}{z-\z_k^j}\right)+A(\xi),
\end{eqnarray}
where
$$
|A(\xi)|\ole\left(\frac{\mb{diam}(I_k^j)}{d(z,I_k^j)}\right)^2.
$$
Since by the definition of $\z_k^j$,
$$
\int_{I_k^j}(\xi-\z_k^j)d\mu_s(\xi)=0,
$$
according to (\ref{2.9})-(\ref{2.12}) we have
\begin{eqnarray}
|\log|P_n(z)|-n\log\kap(K_s)|&\ole&1+
\sum_{j=1}^m\sum_{k=1}^{n_j}\frac{1}{\mu_s(I_k^j)}\int_{I_k^j}
\left(\frac{\mb{diam}(I_k^j)}{d(z,I_k^j)}\right)^2d\mu_s(\xi)
\nonumber\\
\label{2.13}
&=&1+\sum_{j=1}^m\sum_{k=1}^{n_j}
\left(\frac{\mb{diam}(I_k^j)}{d(z,I_k^j)}\right)^2
=:1+\Sigma_3(z).
\end{eqnarray}
Next, we formulate a statement
 which is proved in Section 3.
Let
$$
\Om_{s}:=\{\z\in\Om:s\le g(\z)\le 2s\},\quad s>0.
$$
\begin{lem}\label{lem2.2}
For $z\in  K$ and $n>n_0$ we have
\beq\label{2.14}
\Sigma_3(z)\ole\int_{K_{c/n}}\frac{d(\z,K)}{|\z-z|^2}|d\z|,
\eeq
\beq\label{2.15}
\Sigma_3(z)\ole \int_{\Om_{c/n}}\frac{dm_2(\z)}{|\z-z|^2}.
\eeq
\end{lem}
Thus,  (\ref{2.1}), (\ref{2.13}), and Lemma \ref{lem2.2}
imply the following result which is of independent interest.
\begin{th}\label{th3}
There exist constants $c,c_1,c_2$, and $n_0=n_0(K,c)$ such that for
 $n>n_0$
 we
have
\beq\label{2.16}
||T_n||_K\le c_1{\em{\kap}}(K)^n\sup_{z\in\partial K}
\int_{K_{c/n}}\frac{d(\z,K)}{|\z-z|^2}|d\z|,
\eeq
\beq\label{2.161}
||T_n||_K\le c_2{\em{\kap}}(K)^n\sup_{z\in\partial K}
 \int_{\Om_{c/n}}\frac{dm_2(\z)}{|\z-z|^2}.
\eeq
\end{th}

\absatz{Distortion Properties of  $\phi_j$}

By \cite[p. 23, Lemma 2.3]{andbla}, which is an immediate consequence of
Koebe's one-quarter theorem, we have the following statement.
 Recall that $\psi_j:=\phi_j^{-1}$ is defined in
$
 A^j:=\{w:1<|w|<e^{2s_0/\om_j}\}
 $ and
 $E^j:=\{w:1<|w|<e^{s_0/\om_j}\}$.
\begin{lem}\label{lem3.1}
For $w\in E^j$ and $z=\psi_j(w)$,
$$
c_1^{-1}\frac{d(z,K^j)}{|w|-1}\le|\psi_j'(w)|\le
c_1\frac{d(z,K^j)}{|w|-1}.
$$
Moreover, if $|\tau-w|\le(|w|-1)/2$ and $\z=\psi_j(\tau)$, then
\beq\label{3.2}
c_2^{-1}\frac{|\tau-w|)}{|w|-1}\le\frac{|\z-z|}{d(z,K^j)}\le
c_2\frac{|\tau-w|}{|w|-1}.
\eeq
\end{lem}
In addition to the
 mapping $\phi_j$,
 we consider a conformal (and univalent) mapping $\Phi_j:\C\setminus K^j
 \to\D^*$ normalized by the condition $\Phi_j(\infty)=\infty$
 Then $h_j:=\Phi_j\circ\psi_j$ is a conformal mapping of
 $
 A^j
 $
onto a doubly connected domain bounded by a unit circle and the curve
$h_j(C(e^{2s_0/\om_j}))\subset\D^*$. According to the Carath\'eodory
prime end theorem (see \cite[p. 30, Theorem 2.15]{pom}),
 $h_j$ can be  extended
continuously to the unit circle $C(1)$. Moreover, by the Schwarz reflection
principle (see \cite[p. 4]{pom})
$h_j$  can be extended analytically  into
$\{w: e^{-2s_0/\om_j}<|w|<e^{2s_0/\om_j}\}$. This implies  that
\beq\label{3.3}
\left|\frac{h_j(w_2)-h_j(w_1)}{w_2-w_1}\right|\asymp 1,\quad w_1,w_2\in E^j.
\eeq
Since $\psi_j=\Phi_j^{-1}\circ h_j$ in $E^j$, according to (\ref{3.3}), many known
distortion properties of $\Phi_j$ and $\Phi_j^{-1}$ imply the analogous properties of
$\phi_j$ and $\psi_j$. We describe some of them.

Let for $0<\de\le \de_j:=e^{s_0/(2\om_j)}-1$,
\beq\label{3.1n}
L_\de^j:=K_{\om_j\log(1+\de)}=\{z:|\phi_j(z)|=1+\de\}.
\eeq
For the points $z_k\in L_\de^j$ and $w_k:=\phi_j(z_k)=(1+\de)e^{i\theta_k},
k=1,2$, such that $0<\theta_2-\theta_1<2 \pi$, consider
$$
\hat{L}_\de^j(z_1,z_2):=\{ z=\psi_j((1+\de)e^{i\theta}):\theta_1\le \theta\le\theta_2\},
$$
$$
B_\de^j(z_1,z_2):=\{ z=\psi_j(re^{i\theta}):\theta_1\le \theta\le\theta_2,
0\le r-1-\de\le \min(\theta_2-\theta_1, e^{s_0/\om_j}-1-\de)\}.
$$
By virtue of Lemma \ref{lem3.1},  for $z_1,z_2\in  L^j_\de$ and $\z\in
B_\de^j(z_1,z_2)$
such that  $\theta_2-\theta_1\le \de/c_3$, where $c_3=10c_2(1+\max_j\de_j)$,
we have
$$
\theta_2-\theta_1<\de\le e^{s_0/(2\om_j)}-1<e^{s_0/ \om_j}-e^{s_0/(2\om_j)}
\le e^{s_0/\om_j}-1-\de,
$$
\beq\label{3.3n}
\frac{|\z-z_1|}{d(z_1,K^j)}\le c_2\frac{2(\theta_2-\theta_1)(1+\de_j)}{\de}\le\frac{1}{5},
\eeq
i.e., $|\psi_j'(\z)|\asymp d(z_1,K^j)/\de$
and, therefore,
\begin{eqnarray}
\frac{\theta_2-\theta_1}{c_4\de}d(z_1,K^j)&\le&|z_2-z_1|\le
|\hat{L}_\de^j(z_1,z_2)|\nonumber\\
\label{3.4}
&=&(1+\de)\int_{\theta_1}^{\theta_2}|\psi_j'((1+\de)e^{i\theta})|d\theta
\le
c_4 \frac{\theta_2-\theta_1}{\de}d(z_1,K^j),
\end{eqnarray}
as well as
\begin{eqnarray}
c_5^{-1}\left(\frac{\theta_2-\theta_1}{\de}\right)^2d(z_1,K^j)^2&\le&
m_2(B_\de^j(z_1,z_2))
=\int_{\theta_1}^{\theta_2}\int_{1+\de}^{1+\de+\theta_2-\theta_1}
|\psi_j'(re^{i\theta})|^2rdrd\theta
\nonumber\\
\label{3.5}
&\le&
c_5\left(\frac{\theta_2-\theta_1}{\de}\right)^2d(z_1,K^j)^2.
\end{eqnarray}
{\bf Proof of Lemma \ref{lem2.1}.}
First, we note that by (\ref{3.1n}) $K_s^j=L_\de^j$, where
\beq\label{3.2n}
\frac{s}{\om_j}\le\de=e^{s/\om_j}-1\le e^{s_0/(2\om_j)}\frac{s}{\om_j},\quad
0<\de<\de_j.
\eeq
By (\ref{2.1n}), (\ref{2.5}), and (\ref{3.4}) if $2\pi/n_j\le\de/c_3$
 and $n>n_1=10/(\min_j\om_j)$, then
$$
|I_k^j|\le c_4\frac{2\pi}{n_j\de}d(\xi_{k-1}^j,K^j)
\le\frac{4\pi c_4}{n s}d(\xi_{k-1}^j,K).
$$
Therefore, if we let
$c=100\pi c_4$, then, for $n> n_0:=n_1+n_2+2c/s_0$,
where $n_2$ is defined in (\ref{2.d}),
 and $s=c/n$ we obtain the last inequality in
(\ref{2.7}).
Moreover, (\ref{2.1n}) and the left-hand side of (\ref{3.2}) with $z=\xi_{k-1}^j$
and $\z=\xi_k^j$ yield the first inequality in
(\ref{2.7}).

\newpage

Furthermore, (\ref{2.7}), (\ref{3.3n}), (\ref{3.5}), and (\ref{3.2n})
imply first two order equivalences in (\ref{2.8}).
The last inequality in (\ref{2.8}) follows from (\ref{2.1n}) and (\ref{3.3n})
with $z_1=\xi_{k-1}^j$.

 \hfill$\Box$

{\bf Proof of Lemma \ref{lem2.2}.}
 By (\ref{2.7})
for $z\in  K$ we have
$$
\sum_{k=1}^{n_j}\left(\frac{\mb{diam}(I_k^j)}{d(z,I_k^j)}\right)^2
\ole\sum_{k=1}^{n_j}\frac{|I_k^j|d(I_k^j,K^j)}{d(z,I_k^j)^2}\ole
\int_{K_s^j}\frac{d(\z,K)}{|\z-z|^2}|d\z|,
$$
which implies (\ref{2.14}).

Furthermore, by Lemma \ref{lem2.1} for
$z\in  K$ we obtain
$$
\sum_{k=1}^{n_j}\left(\frac{\mb{diam}(I_k^j)}{d(z,I_k^j)}\right)^2
\ole\sum_{k=1}^{n_j}\frac{m_2(B_k^j)}{d(z,B_k^j)^2}
\ole
\int_{\Om_{s}^j}\frac{dm_2(\z)}{|\z-z|^2},
$$
where $\Om_s^j:=\{\z\in\mb{int}(K_{2s_0}^j):s\le g(\z)\le 2s\}$,
which yields (\ref{2.15}).

 \hfill$\Box$

{\bf Proof of Theorem \ref{th1}.}
Without loss of generality we can assume that $n$ is sufficiently large.
By virtue of (\ref{3.3}) and the Loewner inequality (see \cite{loe}
or \cite[p. 27, Lemma 2.5]{andbla}) for $0<\de<\de_j$  we obtain
$
d(K^j,L_\de^j)\oge \de^{2},
$
which, together with (\ref{3.2n}), imply
\beq\label{3.10}
d_s:=d(K,K_s)\oge s^{2},\quad 0<s<\frac{s_0}{2}.
\eeq
Therefore, for $z\in\partial K$,
$$
\int_{\Om_{s}}\frac{dm_2(\z)}{|\z-z|^2}\le
\int_0^{2\pi}\int_{d_s}^{c_6}\frac{rdrd\theta}{r^2}=
2\pi\log\frac{c_6}{d_s}
$$
and (\ref{1.1}) follows from (\ref{2.161}) and (\ref{3.10}).

\hfill$\Box$

\absatz{Quasismooth Components}

We  need the concept of quasiconformality which can be found in \cite{ahl}
or \cite{lehvir}.
Moreover, almost all facts, necessary for our consideration,
can be derived from a simple statement on the change of the
relative positions of
three points  under a quasiconformal mapping.
We formulate it in a slightly more general form than, for
example, in \cite{andbeldzj, andbla}. However, it is
clear how the proof from there has to be modified
to prove the following result (cf. \cite{abd}).
\begin{lem}\label{lem4.1}(\cite[p. 29, Theorem 2.7]{andbla}).
Let $w=F(\z)$ be a $Q$-quasiconformal mapping of a domain $G\subset\C$
and let $\z_k\in E\subset G,w_k:=F(\z_k),k=1,2,3$, where $E$ is a compact set.
Then, the following is true.

(i) The conditions $|\z_1-\z_2|\le c_1|\z_1-\z_3|$ and $|w_1-w_2|\le
 c_2|w_1-w_3|$ are equivalent. Besides, the constants $c_1$ and $c_2$ are
 mutually dependent and depend on $Q,G,$ and $E$.

 (ii) If $|\z_1-\z_2|\le c_1|\z_1-\z_3|$, then
 $$
 c_3^{-1}\left|\frac{w_1-w_3}{w_1-w_2}\right|^{1/Q}\le
 \left|\frac{\z_1-\z_3}{\z_1-\z_2}\right|\le c_3
\left|\frac{w_1-w_3}{w_1-w_2}\right|^Q,
$$
where $c_3\ge 1$ depends on $c_1,Q,G$, and $E$.
\end{lem}
First, we consider the case where a component $K^j$ of $K$ is a domain bounded by a
quasismooth curve $L^j:=\partial K^j$. By (\ref{1.3}) and Ahlfors' theorem (see
\cite[p. 100, Theorem 8.6]{lehvir}) $L^j$ is quasiconformal and $\phi_j$  can be extended
to a quasiconformal homeomorphism of a neighborhood $N^j=:G$ of $
\psi_j(\ov{E^j})=:E$. Here $E^j=\{ w:1<|w|<e^{s_0/\om_j}\}$.
For $\z\in E$ denote by $\z_0\in L^j$  any point with $|\z-\z_0|=d(\z,L^j)$
and let $\z^*:=\psi_j(\phi_j(\z)/|\phi_j(\z)|)$.
Furthermore, for  $0<\de\le\de_j=e^{s_0/(2\om_j)}-1$
and $z\in L^j$ denote by $z_\de\in L_\de^j$ any point with
$|z-z_\de|=d(z,L^j_\de)$ and let $z_\de^*:=\psi_j((1+\de)\phi_j(z))$.
Applying Lemma \ref{lem4.1}(i) to the triplets $z,z_\de^*,z_\de$ and $\z,\z^*,\z_0$
we obtain
$$
d(z,L^j_\de)=|z-z_\de|\asymp|z-z^*_\de|,
$$
$$
d(\z,L^j)=|\z-\z_0|\asymp|\z-\z^*|.
$$
Moreover, repeated application of Lemma \ref{lem4.1}  implies
that for $z\in L^j$ and $\z\in L^j_\de$,
$$
|\z-z|\le|\z-\z^*|+|\z^*-z|\ole|z_\de^*-\z^*|
$$
and
\begin{eqnarray}
\frac{d(\z,L^j)}{|\z-z|}&\le& \left|\frac{\z-\z^*}{\z-z}\right|
\ole\left|\frac{\phi_j(\z)-\phi_j(\z^*)}{\phi_j(\z)-\phi_j(z)}\right|^{1/Q}
=\left|\frac{\phi_j(z^*_\de)-\phi_j(z)}{\phi_j(z^*_\de)-\phi_j(\z^*)}
\right|^{1/Q}
\nonumber\\
\label{4.4}
&\ole&
\left|\frac{z_\de^*-z}{z_\de^*-\z^*}\right|^{1/Q^2}
\ole
\left(\frac{d(z,L^j_\de)}{|\z-z|}\right)^{1/Q^2}.
\end{eqnarray}
We claim that $L^j_\de$ are ``uniformly" quasismooth. That is,
for $\z_1,\z_2\in L^j_\de$, denote by $L_\de^j(\z_1,\z_2)$ the shorter arc
of $L_\de^j$ between $\z_1$ and $\z_2$.
\begin{lem}\label{lem4.2}
For $\z_1,\z_2\in L^j_\de$ and $0<\de<\de_j$,
\beq\label{4.5}
|L_\de^j(\z_1,\z_2)|\ole |\z_2-\z_1|.
\eeq
\end{lem}
{\bf Proof.} Let
$$\ga:=L_\de^j(\z_1,\z_2)=\{ \psi_j((1+\de)e^{i\theta}):
\theta_1\le\theta\le\theta_2\},\quad 0<\theta_2-\theta_1<2\pi,
$$
$$
\phi_j(\z_k)=:\tau_k=(1+\de)e^{i\theta_k},\quad
k=1,2.
$$
According to (\ref{3.2}) and (\ref{3.4}) there exists
a sufficiently large constant
$c_4$ such that if
$\theta_2-\theta_1\le\de/c_4=:\ve\de$ then
\beq\label{4.6}
|\ga|\ole d(\z_1,L^j)\frac{|\tau_2-\tau_1|}{\de}\asymp|\z_2-\z_1|.
\eeq
Therefore, we can assume that $\theta_2-\theta_1>\ve\de$.
Let $\ga^*:=\{\psi_j(e^{i\theta}): \theta_1\le\theta\le\theta_2\}$.
Consider the points $\theta_1=:\eta_1<\eta_2<\ldots<\eta_p:=\theta_2$
such that
$$
\frac{\ve\de}{2}\le\eta_{k+1}-\eta_k\le\ve\de,\quad k=1,\ldots,p-1.
$$
Let
 $$
 \xi_k:=\psi_j((1+\de)e^{i\eta_k}),\quad \xi_k^*:=\psi_j(e^{i\eta_k}),
 $$
$$
\ga_k:=\{\z=\psi_j((1+\de)e^{i\theta}):\eta_k\le\theta\le\eta_{k+1}\},
$$
$$
\ga_k^*:=\{z=\psi_j(e^{i\theta}):\eta_k\le\theta\le\eta_{k+1}\},
$$
By virtue of Lemma \ref{lem4.1}(i), (\ref{4.6}), and our assumption that $L^j$
is quasismooth we have
$$
|\z_2-\z_1|\asymp|\z_2-\xi_1^*|\asymp|\xi_p^*-\xi_1^*|\asymp|\ga^*|
$$
as well as
$$
|\ga_k|\asymp|\xi_{k+1}-\xi_k|\asymp |\xi_k^*-\xi_k|\asymp
|\xi_k^*-\xi_{k+1}^*|\asymp|\ga_k^*|.
$$
Hence,
$$
|\ga|=\sum_{k=1}^{p-1}|\ga_k|\asymp \sum_{k=1}^{p-1}|\ga_k^*|=|\ga^*|
\asymp|\z_2-\z_1|.
$$

\hfill$\Box$

Further, we claim that
for $z\in\partial K$ and $0<\de\le\de_j$,
\beq\label{4.8}
\int_{L_\de^j}\frac{d(\z,L^j)}{|\z-z|^{2}}|d\z|\ole 1.
\eeq
Indeed, the only
nontrivial  case is where $z\in L^j$.
In this case,  according to Lemma \ref{lem4.2} and
\cite[(3.20)]{and12}, for $\al>0$ we have
\beq\label{4.7}
\int_{L_\de^j}\frac{|d\z|}{|\z-z|^{1+\al}}\le c_5d(z,L_\de^j)^{-\al},
\eeq
where $c_5$ depends on $K$ and $\al$. Furthermore,
by virtue of (\ref{4.4}) and (\ref{4.7}),
$$
\int_{L_\de^j}\frac{d(\z,L^j)}{|\z-z|^{2}}|d\z|\ole
\int_{L_\de^j}\frac{d(z,L^j_\de)^{1/Q^2}}{|\z-z|^{1+1/Q^2}}|d\z|\ole
1
$$
which yields (\ref{4.8}).

Next, let $K^j=:L^j$ be a quasismooth arc.  Below, we show how to modify
the above reasoning   to obtain  (\ref{4.8}) as well.
We only give the main ideas of the proof.
Denote by $z_1$ and $z_2$ the endpoints of $L^j$. For $k=1,2$, let $w_k:=\phi_j(z_k)$,
$$
E^j_1:=\{w:1<|w|<e^{s_0/\om_j},\arg w_1<\arg w<\arg w_2\},
$$
 $$
E^j_2:=E^j\setminus\ov{E_1^j},\quad
R_k^j:=\psi_j(E^j_k),\quad L^j_{\de,k}:= L_\de^j\cap\ov{R_k^j}.
$$
According to the Ahlfors criterion (see \cite[p. 100, Theorem 8.6]{lehvir}) $\partial E^j_k$ is
a quasiconformal curve. Moreover, repeating  the
proof of \cite[p. 30, Lemma 2.8]{andbla} practically word for word, we can show that $\partial
R^j_k$ is also a quasiconformal curve. Therefore,
 $\psi_j$ can be extended quasiconformally into a
neighborhood $N_k^j$
of $\ov{E_k^j}$. Denote this extension by $\psi_{j,k}$.
 Next, we can use Lemma \ref{lem4.1} with $F=\psi_{j,k},
 E=\ov{E^j_k}$, and $G=N_k^j$ to prove the analogues of
 (\ref{4.4}), (\ref{4.5}),
and (\ref{4.7})
with $L_{\de,k}^j$ instead of $L_\de^j$
(for more detail, see the proof of \cite[Lemma 2.3]{andmai}).
This implies (\ref{4.8})  in the case of a quasismooth arc $L^j$
as well.
Summarizing, we have the following result.
\begin{lem}\label{lem4.3}
If the components of $K$ are either closed Jordan domains
bounded by a quasismooth curve or quasismooth arcs, then for
$z\in \partial K$ and $0<s\le s_0/2$,
$$
\int_{K_s}\frac{d(\z,K)}{|\z-z|^2}|d\z|\ole 1.
$$
\end{lem}
{\bf Proof of Theorem \ref{th2}.} The inequality (\ref{1.4})
follows immediately from (\ref{2.16}) and Lemma \ref{lem4.3}.

\hfill$\Box$

  \absatz{Acknowledgements}

 The
author would like to warmly thank    M. Nesterenko for many useful
remarks.

V. V. Andrievskii

 Department of Mathematical Sciences

 Kent State University

 Kent, OH 44242

 USA

e-mail: andriyev@math.kent.edu

tel: 330-672-9029


\begin{thebibliography}{99}


\bibitem{abd}
F. G. Abdullaev (1986): On Orthogonal Polynomials in Domains with Quasiconformal
Boundary. Dissertation, Donetsk (Russian).

\bibitem{ahl}
L. V. Ahlfors (1966): Lectures on Quasiconformal Mappings.
Princeton, N.J.: Van Nostrand.

\bibitem{and12}
V. V. Andrievskii (2012): {\it Weighted $L_p$ Bernstein-type inequalities
on a quasismooth curve in the complex plane.} Acta Math. Hungar., {\bf
135(1-2)}:8--23.

\bibitem{andbeldzj}
V. V. Andrievskii, V. I. Belyi, V. K. Dzjadyk  (1995): Conformal
Invariants in Constructive Theory of Functions of Complex
Variable.  Atlanta, Georgia: World Federation Publisher.


\bibitem{andbla}
  V. V. Andrievskii, H.- P. Blatt (2002): Discrepancy of Signed Measures and
      Polynomial Approximation.
      Berlin/New York: Springer-Verlag.


\bibitem{andmai}
V. V. Andrievskii, V. V. Maimeskul (1995): {\it
Constructive description of certain classes of functions on
quasismooth arcs}. Russian Acad. Sci. Izv. Math., {\bf 44}:193--206.


\bibitem{kovpom}
T. K\"ovary, Ch. Pommerenke (1967): {\it On Faber polynomials and
Faber expansions.} Mathem. Zeitschr., {\bf 99}:193--206.


\bibitem{lehvir}
O. Lehto, K. I. Virtanen (1973):  Quasiconformal Mappings in the Plane, 2nd
ed.,   New York: Springer-Verlag.



 \bibitem{loe}
 K. L\"owner (1919): {\it \"Uber Extremums\"atze der konformen Abbildung des \"Ausseren
 des Einheitskreises.} Math. Z., {\bf 3}:65--77.

\bibitem{lesvinwar}
F. D. Lesley, V. S. Vinge, S. E. Warschawski (1974):
{\it Approximation by polynomials for a class of Jordan
domains.} Mathem. Zeitschr., {\bf 138}:225--237.


\bibitem{pom}
Ch. Pommerenke  (1992):  Boundary Behaviour of Conformal Maps.
Berlin/New York: Springer-Verlag.

\bibitem{ran}
T. Ransford (1995): Potential Theory in the Complex Plane,
Cambridge: Cambridge University Press.


\bibitem{saftot}
 E. B. Saff ,  V. Totik (1997):    Logarithmic Potentials with External
Fields,  New York/Berlin: Springer-Verlag.


\bibitem{smileb}
V. I. Smirnov, N. A. Lebedev (1968):   Functions of a Complex
Variable. Constructive Theory,  Cambridge: Mass. Institute of
Technology.

\bibitem{sodyud}
M. L. Sodin, P. M. Yuditskii (1993): {\it Functions least deviating
from zero on closed subsets of the real line.} St. Petersburg Math. J.,
{\bf 4}:201--249.

\bibitem{sue}
P. K. Suetin (1998): Series of Faber Polynomials, Amsterdam:
Gordon and Breach Science Publishers.

\bibitem{tot12}
V. Totik (2012): {\it Chebyshev polynomials on a system of curves.}
Journal D'Analyse Math\'ematique, {\bf 118}:317--338.

\bibitem{tot13}
V. Totik (2013): {\it Chebyshev polynomials on compact sets.}
Potential Anal., http://dx.doi.org/10.1007/s11118-013-9357-6.


\bibitem{tot14}
V. Totik (2014): {\it Asymptotics of Christoffel functions on
arcs and curves.}
Advances in Mathematics, {\bf 252}:114--149.

\bibitem{wal}
J. L. Walsh (1969):  Interpolation and Approximation by Rational
 Functions in the Complex Plane, 5th ed. Providence, American
 Mathematical Society.



\bibitem{wid}
H. Widom (1969): {\it Extremal polynomials assosiated with a system
of curves in the complex plane.} Adv. Math., {\bf 3}:127--232.




\end{thebibliography}
\end{document}